\documentclass[letterpaper]{article}
\usepackage[english]{babel}

\usepackage{authblk}

\usepackage{amsmath,amsthm}
\usepackage{amssymb}
\usepackage{float}
\usepackage{enumitem}
\usepackage{empheq}

\usepackage{hyperref}

\makeatletter

\newskip\@bigflushglue \@bigflushglue = -100pt plus 1fil

\def\bigcentering{\let\\\@centercr\rightskip\@bigflushglue%
\leftskip\@bigflushglue
\parindent\z@\parfillskip\z@skip}

\makeatother

\usepackage{graphics,graphicx,xcolor}
\usepackage{standalone}
\usepackage{tikz}
\usetikzlibrary{calc,shapes.arrows,decorations.pathmorphing,decorations.pathreplacing,decorations.markings,arrows,positioning,patterns,shapes.multipart}

\definecolor{rouge}{RGB}{255,77,77}
\definecolor{vert}{RGB}{0,178,102}
\definecolor{jaune}{RGB}{255,255,0}
\definecolor{violet}{RGB}{208,32,144}
\definecolor{orange}{RGB}{255,140,0}
\definecolor{bleu}{RGB}{0,0,205}

\theoremstyle{plain}

\newtheorem*{theorem*}{Theorem}
\newtheorem{theorem}{Theorem}

\newtheorem{proposition}[theorem]{Proposition}

\newtheorem*{questions*}{Questions}

\newtheorem{remark}[theorem]{Remark}

\def\BB{\overrightarrow{B}}
\def\H{\mathbb{H}}
\def\N{\mathbb{N}}

\def\R{\mathbb{R}}
\def\Z{\mathbb{Z}}
\def\bb{\overrightarrow{b}}
\def\xx{\overrightarrow{x}}
\def\yy{\overrightarrow{y}}
\def\zz{\overrightarrow{z}}
\def\ellell{\overrightarrow{\ell}}
\def\rr{\overrightarrow{r}}

\newcommand{\define}[1]{\emph{#1}}

\title{Addendum to "\emph{Tilings problems on Baumslag-Solitar groups}"}
\date{}
\author[1]{Nathalie Aubrun}
\author[2]{Jarkko Kari}
\affil[1]{CNRS, Université Paris-Saclay, LISN, 91400 Orsay, France}
\affil[2]{Department of Mathematics, University of Turku, FIN-20014, Turku, Finland}
\begin{document}
 
\maketitle 
 
\begin{abstract}
In our article~\cite{AubrunKari2013} we state the the Domino problem is undecidable for all Baumslag-Solitar groups $BS(m,n)$, and claim that the proof is a direct adaptation of the construction of a weakly aperiodic subshift of finite type for $BS(m,n)$ given in the paper. In this addendum, we clarify this point and give a detailed proof of the undecidability result.  We assume the reader is already familiar with the article~\cite{AubrunKari2013}.
\end{abstract}

\section*{Introduction}
\label{section.introduction}

In~\cite{AubrunKari2013} we state as a direct corollary of the main construction that the Domino problem is undecidable on all Baumslag-Solitar groups $BS(m,n)$. It turns out that it is not as immediate as we write it, and we believe that this result deserves a full explanation.

The proof is based on the proof of the undecidability of the Domino problem on the discrete hyperbolic plane given by the second author in~\cite{Kari2007}.This latter is an adaptation of a former construction of a strongly aperiodic SFT on~$\Z^2$~\cite{Kari1996}. This proof proceeds by reduction to the immortality problem for rational piecewise affine maps. We first recall the key ingredients of this proof.

    A mapping $f:U\to U\subset \R^2$ is a rational piecewise affine map if there exists a partition $U=U_1\cup U_2 \cup\dots\cup U_n$ where every  $U_i$ is a unitary square with integer coordinates, and such that $f=f_i$ on every $U_i$,  and $f_i:U_i\to\R^2$ is an affine function with rational parameters. A point $\xx\in U$ is immortal for $f$ if for every $k\in\Z$, the iterated image $f^k(\xx)$ belongs to $U$. The immortality problem for rational piecewise affine maps is the decision problem that inputs such a function $f$, and outputs $\texttt{Yes}$ if $f$ possesses an immortal point, and $\texttt{No}$ otherwise. This problem reduces to the immortality problem for Turing machines, which is known to be undecidable~\cite{Hooper1966}.

    \begin{theorem}[\cite{Kari2007}]
    The immortality problem for rational piecewise affine maps is undecidable.
    \end{theorem}

    The second author proves in~\cite{Kari2007} that the problem of tiling the discrete hyperbolic plane $\H_2$ with pentagonal Wang tiles is undecidable, by a reduction to the immortality problem for rational piecewise affine maps. The proof is based on the construction, for every piecewise affine map $f:U\to U$ with rational parameters, of a finite tileset that computes the function $f$, meaning that a tiling by this tileset encodes the orbit of a point $\xx\in U$ under the action of~$f$.

    \begin{theorem}[\cite{Kari2007}]
    \label{DP_indecidable_hyperbolique}
    The Domino problem is undecidable on the discrete hyperbolic plane $\H_2$.
    \end{theorem}

    The most difficult part in the construction is to ensure finiteness of the tileset. This issue may be bypassed by combining two main ingredients: representing real numbers by Beatty sequences, and taking advantage of the rationality of function $f$. Instead of Beatty sequences of a point $\xx\in\R^2$ we use its balanced representation, to be defined on page~\pageref{definition:representation_equilibree}, and that takes accounts of the merges between sheets in $BS(m,n)$.

    \section*{Adaptation to Baumslag-Solitar groups \texorpdfstring{$BS(m,n)$}{BSmn}}
    \label{section.adaptation_BSmn}

    Following~\cite{Kari2007} we prove that the Domino problem is undecidable for all Baumslag-Solitar $BS(m,n)$ for all integers $m,n\in\N^*$

    \[BS(m,n)=\langle a,t | t^{-1}a^mt=a^n \rangle.\]
    
Since $BS(-m,-n)$ is isomorphic to $BS(m,n)$, it is enough to consider groups with $m>0$. For simplicity, we also assume that $n>0$. The case $n<0$ is analogous. The Cayley graph of $BS(m,n)$ with generating set $\{a,t,a^{-1},t^{-1}\}$ is made of several sheets that merge $m$ by $m$ from top to give $n$ other sheets to the bottom, so that the global structure these sheets are arranged looks like an $(m+n)$-regular tree. Each of these sheets is quasi-isometric to the hyperbolic plane $\mathbb{H}^2$.    

    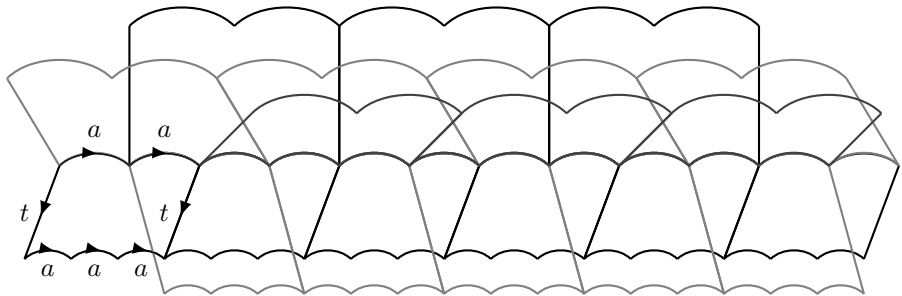
\begin{figure}[!ht]
        \centering
\begin{tikzpicture}[scale=0.31]
\foreach \x in {0,...,5} {
\draw[thick] (6*\x-1.5,0) -- (6*\x,4);
\draw[thick] (6*\x,4) to [controls=+(45:1) and +(135:1)] (6*\x+3,4);
\draw[thick] (6*\x+3,4) to [controls=+(45:1) and +(135:1)] (6*\x+6,4);
\draw[thick] (6*\x-1.5,0) to [controls=+(45:0.7) and +(135:0.7)] (6*\x+2-1.5,0);
\draw[thick] (6*\x+2-1.5,0) to [controls=+(45:0.7) and +(135:0.7)] (6*\x+4-1.5,0);
\draw[thick] (6*\x+4-1.5,0) to [controls=+(45:0.7) and +(135:0.7)] (6*\x+6-1.5,0);
\draw[thick] (6*\x+6,4) -- (6*\x+6-1.5,0);
}

\foreach \x in {0,...,4} {
\begin{scope}[shift={(3,0)},color=black!50]
\draw[thick] (6*\x+1.5,0-1.5) -- (6*\x,4);
\draw[thick] (6*\x,4) to [controls=+(45:1) and +(135:1)] (6*\x+3,4);
\draw[thick] (6*\x+3,4) to [controls=+(45:1) and +(135:1)] (6*\x+6,4);
\draw[thick] (6*\x+1.5,0-1.5) to [controls=+(45:0.7) and +(135:0.7)] (6*\x+2+1.5,0-1.5);
\draw[thick] (6*\x+2+1.5,0-1.5) to [controls=+(45:0.7) and +(135:0.7)] (6*\x+4+1.5,0-1.5);
\draw[thick] (6*\x+4+1.5,0-1.5) to [controls=+(45:0.7) and +(135:0.7)] (6*\x+6+1.5,0-1.5);
\draw[thick] (6*\x+6,4) -- (6*\x+6+1.5,0-1.5);
\end{scope}
}

\foreach \x in {0,1,2} {
\begin{scope}[scale=3/2,shift={(2,4*2/3)}]
\draw[thick] (6*\x,0) -- (6*\x,4);
\draw[thick] (6*\x,4) to [controls=+(45:1) and +(135:1)] (6*\x+3,4);
\draw[thick] (6*\x+3,4) to [controls=+(45:1) and +(135:1)] (6*\x+6,4);
\draw[thick] (6*\x,0) to [controls=+(45:0.7) and +(135:0.7)] (6*\x+2,0);
\draw[thick] (6*\x+2,0) to [controls=+(45:0.7) and +(135:0.7)] (6*\x+4,0);
\draw[thick] (6*\x+4,0) to [controls=+(45:0.7) and +(135:0.7)] (6*\x+6,0);
\draw[thick] (6*\x+6,4) -- (6*\x+6,0);
\end{scope}
}

\foreach \x in {0,1,2,3} {
\begin{scope}[scale=3/2,shift={(0,4*2/3)},color=black!50]
\draw[thick] (6*\x,0) -- (6*\x-1.5,4-1.5);
\draw[thick] (6*\x-1.5,4-1.5) to [controls=+(45:1) and +(135:1)] (6*\x+3-1.5,4-1.5);
\draw[thick] (6*\x+3-1.5,4-1.5) to [controls=+(45:1) and +(135:1)] (6*\x+6-1.5,4-1.5);
\draw[thick] (6*\x,0) to [controls=+(45:0.7) and +(135:0.7)] (6*\x+2,0);
\draw[thick] (6*\x+2,0) to [controls=+(45:0.7) and +(135:0.7)] (6*\x+4,0);
\draw[thick] (6*\x+4,0) to [controls=+(45:0.7) and +(135:0.7)] (6*\x+6,0);
\draw[thick] (6*\x+6-1.5,4-1.5) -- (6*\x+6,0);
\end{scope}
}

\foreach \x in {0,1,2} {
\begin{scope}[scale=3/2,shift={(4,4*2/3)},color=black!75]
\draw[thick] (6*\x,0) -- (6*\x+1.5,4-2.5);
\draw[thick] (6*\x+1.5,4-2.5) to [controls=+(45:1) and +(135:1)] (6*\x+3+1.5,4-2.5);
\draw[thick] (6*\x+3+1.5,4-2.5) to [controls=+(45:1) and +(135:1)] (6*\x+6+1.5,4-2.5);
\draw[thick] (6*\x,0) to [controls=+(45:0.7) and +(135:0.7)] (6*\x+2,0);
\draw[thick] (6*\x+2,0) to [controls=+(45:0.7) and +(135:0.7)] (6*\x+4,0);
\draw[thick] (6*\x+4,0) to [controls=+(45:0.7) and +(135:0.7)] (6*\x+6,0);
\draw[thick] (6*\x+6+1.5,4-2.5) -- (6*\x+6,0);
\end{scope}
}

\begin{scope}[decoration={
    markings,
    mark=at position 0.55 with {\arrow[scale=1.5]{latex}}}
    ]
\draw[postaction={decorate}] (0,4) -- (-1.5,0);
\draw node at (-1.5,2) {$t$};   
\draw[postaction={decorate}] (0,4) to [controls=+(45:1) and +(135:1)] (3,4);
\draw node at (1.5,5.5) {$a$};   
\draw[postaction={decorate}] (3,4) to [controls=+(45:1) and +(135:1)] (6,4);
\draw node at (4.5,5.5) {$a$};  
\draw[postaction={decorate}] (6,4) -- (4.5,0);
\draw node at (4.5,2) {$t$};   
\end{scope}
   
\begin{scope}[decoration={
    markings,
    mark=at position 0.65 with {\arrow[scale=1.5]{latex}}}
    ]   
\draw[postaction={decorate}] (-1.5,0) to [controls=+(45:0.7) and +(135:0.7)] (0.5,0);
\draw[postaction={decorate}] (0.5,0) to [controls=+(45:0.7) and +(135:0.7)] (2.5,0);
\draw[postaction={decorate}] (2.5,0) to [controls=+(45:0.7) and +(135:0.7)] (4.5,0);
\draw node at (-0.5,-0.5) {$a$};  
\draw node at (1.5,-0.5) {$a$};  
\draw node at (3.5,-0.5) {$a$};  
\end{scope}

\end{tikzpicture}
        \caption{A portion of the right Cayley graph of the group $BS(2,3)$ with generating set  $\{a,t,a^{-1},t^{-1}\}$. Three sheets from the top merge and separate into two sheets to the bottom.}
            \label{figure.Cayley_graph_BS23}
    \end{figure}

%


        The group $BS(m,n)$ is embedded into $\R^2$ through a function $\Phi_{m,n}:BS(m,n)\to\mathbb{R}^2$, that is first defined recursively on words $w$ on the alphabet $A=\{ a,t,a^{-1},t^{-1}\}$. If $x$ is a letter from $A$, denote $|w|_x$ the number of occurrences of the letter $x$ in the word $w$. We also call \define{contribution} of $x$ to $w$ the integer $\parallel w\parallel_x=|w|_x-|w|_{x^{-1}}$. With these notations we first define a function $\beta:A^*\to\Z$ by $\beta(w):=-\parallel w\parallel_{t}$, and a function $\alpha_{m,n}:A^*\to\mathbb{R}$, simply denoted $\alpha$ in the sequel, which is defined by induction on the length of words ($\varepsilon$ denotes the empty word) by:
        
        \begin{align*}
        \alpha(\varepsilon) &= 0\\
        \alpha(w.t) &=\alpha(w.t^{-1})=\alpha(w)\\
        \alpha(w.a) &=\alpha(w)+{\left(\frac{m}{n}\right)}^{-\beta(w)}\\
        \alpha(w.a^{-1}) &=\alpha(w)-{\left(\frac{m}{n}\right)}^{-\beta(w)}
        \end{align*}

By induction on the length of words $w$ we get the formula:
        
        \medskip
        
        \begin{proposition}
        \label{proposition.alpha_relation_recurrence}
        For every words $u,v\in A^*$ one has
        \[\alpha(u.v)=\alpha(u)+\left(\frac{m}{n}\right)^{-\beta(u)}\alpha(v).\]
        \end{proposition}
        
        In the sequel we will use in particular the following equalities:
        \begin{align*}
         \alpha(ga) &=\alpha(g)+\left(\frac{m}{n}\right)^{-\beta(g)}  \\
         \beta(gt) &= \beta(g)-1
        \end{align*}

        Finally the function $\Phi_{m,n}:BS(m,n)\to \mathbb{R}^2$ is
        \[\Phi_{m,n}(g)=\left(\alpha(w) ,\beta(w) \right),\]
        where $w$ is a word that represents the group element $g$. One can check that $\Phi_{m,n}$ is well-defined, i.e. the value for $\Phi_{m,n}(g)$ does not depend on the word $w$ chosen, thanks to Proposition~\ref{proposition.alpha_relation_recurrence}.


        \medskip
        
        \begin{proposition}
        The function $\Phi_{m,n}$ is well-defined on $BS(m,n)$.
        \end{proposition}

Note that for $m=1$ we find the same function as the isomorphism $\Phi$ defined only for amenable Baumslag-Solitar groups in~\cite{AubrunSchraudner2020}.       
        
        
        \begin{remark}
        If $|m|\neq 1$ and $|n|\neq 1$ then the function $\Phi_{m,n}$ is not injective. In~\cite{AubrunKari2013} we give an example of injectivity default for $m=3$ and $n=2$: the group element $\omega=bab^{-1}a^2ba^{-1}b^{-1}a^{-2}$ is sent to the origin by$\Phi_{3,2}$, but has infinite order. In~\cite{EsnayMoutot2020} the authors exhibit the word $\omega=bab^{-1}aba^{-1}b^{-1}a^{-1}$ which satisfies that $\Phi_{m,n}(\omega)=(0,0)$ for all $m,n$ such that $|m|\neq 1$ and $|n|\neq 1$.
        \end{remark}

        
        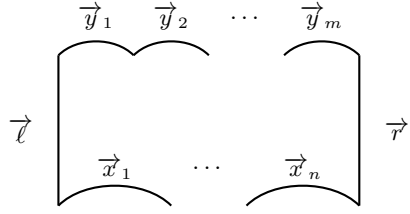
\begin{figure}[!h]
        \centering
\begin{tikzpicture}[scale=0.5]
 
\draw[thick] (0,0) -- (0,4);
\draw[thick] (0,4) to [controls=+(45:0.7) and +(135:0.7)] (2,4);
\draw[thick] (2,4) to [controls=+(45:0.7) and +(135:0.7)] (4,4);
\draw[thick] (6,4) to [controls=+(45:0.7) and +(135:0.7)] (8,4);
\draw[thick] (8,4) -- (8,0);
\draw[thick] (0,0) to [controls=+(45:1) and +(135:1)] (3,0);
\draw[thick] (5,0) to [controls=+(45:1) and +(135:1)] (8,0);


\draw (1.5,1) node{\small$\xx_1$};
\draw (4,1) node{\small$\dots$};
\draw (6.5,1) node{\small$\xx_n$};
\draw (-1,2) node{\small$\ellell$};
\draw (9,2) node{\small$\rr$};
\draw (1,5) node{\small$\yy_1$};
\draw (3,5) node{\small$\yy_2$};
\draw (5,5) node{\small$\dots$};
\draw (7,5) node{\small$\yy_m$};

\end{tikzpicture}
        \caption{A Wang tile for $BS(m,n)$.}
        \label{figure.tuile_BS_calcule}
        \end{figure}
        
        Fix two integers $m,n\in\N^*$. A tile on $BS(m,n)$ \define{computes} a function $f_i:U_i\subset\R^2\to\R^2$ if, the following holds (colors on the edges of the tile are named after Figure~\ref{figure.tuile_BS_calcule}):
        \[\frac{\yy_1+\dots \yy_m}{m}+\rr=f_i\left(\frac{\xx_1+\dots \xx_n}{n}\right)+\ellell.\]
        
        For a single sheet of $BS(m,n)$, one can easily adapt what is done for the discrete hyperbolic plane $\H_2$~\cite{Kari2007}: select all tiles satisfying the relation with colors on the edges belonging to a well chosen finite set. For the whole group difficulties arise where different sheets merge: clearly the direct adaptation of the~$\H_2$ case is not enough, and the tileset should be enriched to take into account the specific structure of $BS(m,n)$. Our solution uses function $\Phi_{m,n}$, from which we define a function $\lambda:BS(m,n)\to\R$ by
        \[\lambda(g):=\frac{1}{m}\left(\frac{n}{m}\right)^{-\beta(g)}\alpha(g),\]
        for every $g\in BS(m,n)$, and one can check that the following holds
        \begin{align*}
         \lambda(ga) &= \lambda(g) + \frac{1}{m}\\
         \lambda(gt) &= \frac{n}{m}\lambda(g).
        \end{align*}

        Thanks to the function $\lambda$ and to the properties it satisfies, we detail the content of every tile that computes a piecewise affine function $f_i:U_i\subset\R^2\to\R^2$ such that $f_i(\xx)=M\xx+\bb$:
        \begin{align*}
         \xx_k(g,\xx) &:= \lfloor \left(n\lambda(g)+k \right)\xx \rfloor - \lfloor \left(n\lambda(g)+(k-1) \right)\xx \rfloor \text{~~for }k=1\dots m\\         
         \yy_k(g,\xx) &:= \lfloor \left(m\lambda(g)+k \right)f_i(\xx) \rfloor - \lfloor \left(m\lambda(g)+(k-1) \right)f_i(\xx) \rfloor \text{~~for }k=1\dots n\\         
         \ellell(g,\xx) &:= \frac{1}{n}f_i\left(\lfloor n\lambda(g)\xx \rfloor\right) - \frac{1}{m}\lfloor m\lambda(g)f_i(\xx) \rfloor + \lfloor \lambda(g)-\frac{1}{2}\rfloor \bb \\
         \rr(g,\xx) &:= \frac{1}{n}f_i\left(\lfloor \left(n\lambda(g)+n\right)\xx \rfloor\right) - \frac{1}{m}\lfloor \left(m\lambda(g)+m\right)f_i(\xx) \rfloor + \lfloor \lambda(g)+\frac{1}{2}\rfloor \bb
        \end{align*}
        
        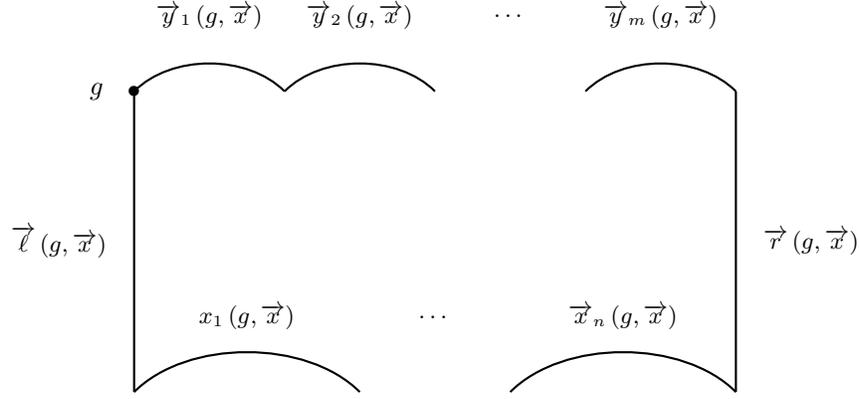
\begin{figure}[!h]
        \centering
\begin{tikzpicture}[scale=1,every node/.style={scale=1}
]
\draw[thick] (0,0) -- (0,4);
\draw[thick] (0,4) to [controls=+(45:0.7) and +(135:0.7)] (2,4);
\draw[thick] (2,4) to [controls=+(45:0.7) and +(135:0.7)] (4,4);
\draw[thick] (6,4) to [controls=+(45:0.7) and +(135:0.7)] (8,4);
\draw[thick] (8,4) -- (8,0);
\draw[thick] (0,0) to [controls=+(45:1) and +(135:1)] (3,0);
\draw[thick] (5,0) to [controls=+(45:1) and +(135:1)] (8,0);


\draw (1.5,1) node{\small$x_1\left(g,\xx\right)$};
\draw (4,1) node{\small$\dots$};
\draw (6.5,1) node{\small$\xx_n\left(g,\xx\right)$};
\draw (-1,2) node{\small$\ellell\left(g,\xx\right)$};
\draw (9,2) node{\small$\rr\left(g,\xx\right)$};
\draw (1,5) node{\small$\yy_1\left(g,\xx\right)$};
\draw (3,5) node{\small$\yy_2\left(g,\xx\right)$};
\draw (5,5) node{\small$\dots$};
\draw (7,5) node{\small$\yy_m\left(g,\xx\right)$};

\draw (0,4) node{$\bullet$};
\draw (-0.5,4) node{$g$};

\end{tikzpicture}
        \caption{Tile 
        to encode a piecewise affine map $f_i:U_i\subset\R^2\to\R^2$ on the group $BS(m,n)$.}
        \label{figure.exemple_bloc_codant}
        \end{figure}
        
        We check that the tile on Figure~\ref{figure.exemple_bloc_codant} does compute the function $f_i$, in other words that the quantity $S$ detailed below sums to null vector $\overrightarrow{0}$:
        \begin{align*}
         S &:= \frac{\yy_1+\dots +\yy_m}{m}+\rr-f_i\left(\frac{\xx_1+\dots +\xx_n}{n}\right)-\ellell
        \end{align*}

        By replacing every term $\yy_k$ and $\xx_k$ by its expression given above, the two sums $\yy_1+\dots +\yy_m$ and $\xx_1+\dots +\xx_n$ telescope and $S$ simplifies in
        
        \begin{align*}
         S &= \frac{1}{m}\lfloor \left(m\lambda(g)+m \right)f_i(\xx) \rfloor - \frac{1}{m}\lfloor m\lambda(g)f_i(\xx) \rfloor + \frac{1}{n}f_i\left(\lfloor \left(n\lambda(g)+n\right)\xx \rfloor\right) - \frac{1}{m}\lfloor \left(m\lambda(g)+m\right)f_i(\xx) \rfloor \\
         & + \lfloor \lambda(g)+\frac{1}{2}\rfloor \bb - f_i\left(\frac{1}{n}\lfloor \left(n\lambda(g)+n \right)\xx \rfloor - \frac{1}{n}\lfloor n\lambda(g)\xx \rfloor\right)\\
          &  - \frac{1}{n}f_i\left(\lfloor n\lambda(g)\xx \rfloor\right) + \frac{1}{m}\lfloor m\lambda(g)f_i(\xx) \rfloor - \lfloor \lambda(g)-\frac{1}{2}\rfloor \bb
        \end{align*}
        
        which then reduces to
        
        \begin{align*}
         S &= \frac{1}{n}f\left(\lfloor \left(n\lambda(g)+n\right)\xx \rfloor\right) + \lfloor \lambda(g)+\frac{1}{2}\rfloor \bb\\
          & - f\left(\frac{1}{n}\lfloor \left(n\lambda(g)+n \right)\xx \rfloor - \frac{1}{n}\lfloor n\lambda(g)\xx \rfloor\right) - \frac{1}{n}f\left(\lfloor n\lambda(g)\xx \rfloor\right) - \lfloor \lambda(g)-\frac{1}{2}\rfloor \bb.
        \end{align*}
        
        We now use the fact that $f_i(c\yy-c\zz)=cf_i(\yy)-cf_i(\zz)+\bb$ to obtain:
        
        \begin{align*}
         S &= - \frac{1}{n}f_i\left(\lfloor n\lambda(g)\xx \rfloor\right) + \lfloor \lambda(g)+\frac{1}{2}\rfloor \bb\\
          & - \frac{1}{n}f_i\left(\lfloor \left(n\lambda(g)+n \right)\xx \rfloor\right) + \frac{1}{n}f_i\left(\lfloor n\lambda(g)\xx \rfloor\right) - \bb + \frac{1}{n}f_i\left(\lfloor \left(n\lambda(g)+n\right)\xx \rfloor\right) - \lfloor \lambda(g)-\frac{1}{2}\rfloor \bb
        \end{align*}
        
        afterwards only remains:
        
        \begin{align*}
         S &= \lfloor \lambda(g)+\frac{1}{2}\rfloor \bb- \bb - \lfloor \lambda(g)-\frac{1}{2}\rfloor \bb
        \end{align*}
        
        and since $\lfloor z+\frac{1}{2} \rfloor - \lfloor z-+\frac{1}{2} \rfloor=1$ for every real number $z$, we finally conclude that $S=0$. 
        
        \medskip
        
        \begin{proposition}
        \label{proposition:tuile_BSmn_calcule_f}
         For every $g\in BS(m,n)$ and every $\xx\in U_i$, the tile described on Figure~\ref{figure.exemple_bloc_codant} computes the piecewise affine map $f_i:U_i\subset \R^2\to\R^2$.
        \end{proposition}
        
        Every tile thus individually computes the image by $f_i$ of the average of the elements on the bottoms edges, and redistributes this image on the top edges. This is performed up to calculation errors, that are stored in the left and right edges of the tile.
        
        As explained in~\cite{AubrunKari2013}, for all $\xx\in\R^2$ and $z\in\R$, if one defines for every $k\in\Z$ 
        \[\BB_k(\xx,z):=\lfloor\left(z+k\right)\xx\rfloor-\lfloor\left(z+(k-1)\right)\xx\rfloor,\]
        then the bi-infinite sequence $\left(\BB_k(\xx,z)\right)_{k\in\Z}$ is a \define{balanced representation of $\xx$}\label{definition:representation_equilibree}. In particular, it is a representation of $\xx=(x_1,x_2)$, meaning that 
        \begin{itemize}
         \item every $\BB_k(\xx,z)$ has integer coordinates in $\left\{\lfloor x_1\rfloor ;\lfloor x_1\rfloor+1\right\}\times\left\{\lfloor x_2\rfloor ;\lfloor x_2\rfloor+1 \right\}$ ;
         \item the following average converges towards $\xx$
        \[\lim_{k\rightarrow\infty} \frac{1}{2k+1}\sum_{j=-k}^{k} \BB_j(\xx,z) =\xx.\]
        \end{itemize}

        \begin{proposition}
        \label{proposition:ligne_tuiles_BSmn_calcule_f}
         For every $g_0\in BS(m,n)$ and every $\xx\in U_i$, if we put the tile from Figure~\ref{figure.exemple_bloc_codant} in position $g$ for every $g\in \left\{ g_0\cdot a^k\mid k\in\Z\right\}$, then one can read the balanced representation $\BB_k\left(\xx,\lambda(g)\right)$ of $\xx$ on the bottom edges and the balanced representation $\BB_k\left(f_i(\xx),\lambda(gt^{-1})\right)$ of $f_i(\xx)$ on the top edges.
        \end{proposition}
        
        Proposition~\ref{proposition:ligne_tuiles_BSmn_calcule_f} expresses the fact that moving from a single tile that computes $f_i$ with errors to an infinite row of tiles makes the calculation of $f_i$ exact.
        
        \medskip
        
        We now check that among all possible tiles that compute $f_i$, we can restrict to a finite tileset. For every $g\in BS(m,n)$, the sequence $\BB_k(\xx,\lambda(g))$ is a balanced representation of $\xx$, so that there exist only finitely many possible values for $\xx_k$, and the same argument prevails for $\yy_k$. It remains to check that the $\ellell$ and $\rr$ can be chosen among a finite set. Using the fact that $\lambda(ga^m)=\lambda(g)+1$, we remark that $\ellell(ga^m,\xx)=\rr(g,\xx)$ ; hence it is enough to ensure a finite number of choices for the $\ellell$ only. Remind that
        
        \[\ellell(g,\xx) := \frac{1}{n}f_i\left(\lfloor n\lambda(g)\xx \rfloor\right) - \frac{1}{m}\lfloor m\lambda(g)f_i(\xx) \rfloor + \lfloor \lambda(g)-\frac{1}{2}\rfloor \bb.\]
        
        We first check that $\ellell(g,\xx)$ is bounded, as a consequence of $\zz-1\leq \lfloor \zz \rfloor <\zz$ (inequalities shall apply coordinate by coordinate). Indeed:

        \begin{align*}
        \lambda(g)\left(f_i(\xx)-\bb\right)-\frac{1}{n}M\overrightarrow{1} -\frac{1}{n}\bb-\lambda(g)f_i(\xx)+\lambda(g)\bb--\frac{1}{2}\bb + \bb&< \ellell(g,\xx) \\
        -\frac{1}{n}M\overrightarrow{1} - \frac{3n+2}{2n}\bb&< \ellell(g,\xx)
        \end{align*}
        
        and

        \begin{align*}
         \ellell(g,\xx) &< \lambda(g)\left(f_i(\xx)-\bb\right)+\frac{1}{n}\bb -\lambda(g)f_i(\xx)+\lambda(g)\bb--\frac{1}{m}\overrightarrow{1}+ \lambda(g)\bb -\frac{1}{2} \bb \\
         \ellell(g,\xx) &< -\frac{1}{m}\overrightarrow{1} - \frac{n-2}{2n}\bb.
        \end{align*}        
        
        Since both vector $\bb$ and matrix $M$ have rational coefficients, one can put at the same denominator $q$ the two inequalities above,  so that there exist two vectors $\overrightarrow{p_1},\overrightarrow{p_2}\in\mathbb{Z}^2$ such that 
        \[
        \frac{\overrightarrow{p_1}}{q} \leq \ellell(g,\xx) \leq \frac{\overrightarrow{p_2}}{q},
        \]
        where $\overrightarrow{p_1}$ is chosen maximal and $\overrightarrow{p_2}$ minimal. Better than that, the value for $\ellell(g,\xx)$ should belong to the finite set
        \[
        \left\{ \frac{\overrightarrow{p_1}}{q},\frac{\overrightarrow{p_1}+(0,1)}{q},\frac{\overrightarrow{p_1}+(1,0)}{q},\frac{\overrightarrow{p_1}+\overrightarrow{1}}{q},\dots,\frac{\overrightarrow{p_2}}{q}\right\} \subset \mathbb{Q}
        \]
        for every $g\in BS(m,n)$ and every $\xx\in U$. Indeed, a careful observation of rational numbers that appear in the expression of $\ellell(g,\xx)$, shows that $\ellell(g,\xx)$ can be written as $\frac{\overrightarrow{p}}{q}$. The fact that $\overrightarrow{p_1}\leq \overrightarrow{p} \leq \overrightarrow{p_2}$ directly follows from the definition of $\overrightarrow{p_1}$ and $\overrightarrow{p_2}$. The tileset $\tau_{f_i}$ corresponding to the function $f_i$ is thus finite.
        
        \medskip
        
        \begin{proposition}
         There exists a finite number of tiles on $BS(m,n)$ with colors as on Figure~\ref{figure.exemple_bloc_codant} that computes $f_i(\xx)$ for every $\xx\in U_i$.
        \end{proposition}
        
        Thanks to the properties of the function $\lambda$ stated above, one has
        \begin{align*}
        y_1\left(g\cdot a^k,\xx\right) &= y_{1+k}\left(g,\xx\right)\text{ for }k\in[1;m-1]\\
        y_1\left(gt,\xx\right) &= x_1\left(g,f_i\left(\xx\right)\right), 
        \end{align*}
        which ensures that for a given $\xx\in U_i$, there exists a tiling of the coset $\{a^k\mid k\in\Z\}$ such that the balanced representations of $\xx$ and $f_i(\xx)$ appear respectively on bottom and top edges. We then put together all tilesets corresponding to every function $f_i$, by adding to these tiles the number $i$ of the function $f_i$ they encode. With the additional local rule that two tiles in positions $g$ and $ga$ should share the same number $i$, we finally get the desired result.
        
        \medskip
        
        \begin{theorem}
		The Domino problem is undecidable on  Baumslag-Solitar groups $BS(m,n)$ for every integers $m,n\in\Z$.
		\end{theorem}

\section*{Acknowledgments}

This work was partially supported by the ANR project CoCoGro (ANR-16-CE40-0005).

\bibliographystyle{alpha}
\bibliography{biblio}

\begin{thebibliography}{Hoo66}

\bibitem[AK13]{AubrunKari2013}
Nathalie Aubrun and Jarkko Kari.
\newblock {Tiling Problems on Baumslag-Solitar groups.}
\newblock In Turlough Neary and Matthew Cook, editors, {\em { Proceedings}
  Machines, Computations and Universality 2013, { Z\"urich, Switzerland,
  9/09/2013 - 11/09/2013}}, volume 128 of {\em Electronic Proceedings in
  Theoretical Computer Science}, pages 35--46. Open Publishing Association,
  2013.

\bibitem[AS20]{AubrunSchraudner2020}
Nathalie Aubrun and Michael Schraudner.
\newblock Tilings of the hyperbolic plane of substitutive origin as subshifts
  of finite type on baumslag-solitar groups {$BS(1,n)$}.
\newblock \url{https://arxiv.org/abs/2012.11037}, 2020.

\bibitem[EM20]{EsnayMoutot2020}
Julien Esnay and Etienne Moutot.
\newblock {Weakly and Strongly Aperiodic Subshifts of Finite Type on
  Baumslag-Solitar Groups}.
\newblock \url{https://arxiv.org/abs/2004.02534}, 2020.

\bibitem[Hoo66]{Hooper1966}
Philip Hooper.
\newblock The undecidability of the {T}uring machine immortality problem.
\newblock {\em The Journal of Symbolic Logic}, 31:219--234, 1966.

\bibitem[Kar96]{Kari1996}
J.~Kari.
\newblock A small aperiodic set of wang tiles.
\newblock {\em Discrete Mathematics}, 160:259 -- 264, 1996.

\bibitem[Kar07]{Kari2007}
Jarkko Kari.
\newblock The tiling problem revisited.
\newblock In {\em Proceedings of the 5th International Conference on Machines,
  Computations, and Universality}, MCU'07, page 72–79, Berlin, Heidelberg,
  2007. Springer-Verlag.

\end{thebibliography}
 
\end{document}